%From carlen@math.gatech.edu  Fri Jan 12 10:38:40 2007
%The change suggested by Carlen was made.
%The reference on the title page to AMS was updated.

\iffalse
Happy New Year!   Here is the final version of the paper. I checked this
against the published version, and it is the same -- apart from formatting
the info in the footnotes. In the course of checking, I did note that we 
should probaly
add a dropped word (not so important) in line -10, page 2 of the output 
the tex file
produces:

not true in general positive operators   ---->   not true in general for 
positive operators

Had we published with Springer, it would be impossible to chenge this 
(without
violating the "definitive version" clauses in their publishing 
agreements. But
here I suppose we can, no?)

\fi

\overfullrule=0pt
\magnification=1200
\baselineskip =1.5\baselineskip
\font\eightit=cmti8
\outer\def\beginsection#1\par{\vskip0pt plus.3\vsize\penalty-100
   \vskip0pt plus-.3\vsize\bigskip\vskip\parskip
   \message{#1}\leftline{\bf#1}\nobreak\smallskip}

\def\IR{{\rm I\kern -1.6pt{\rm R}}}
\def\IP{{\rm I\kern -1.6pt{\rm P}}}
\def\ZZ{{\rm Z\kern -4.0pt{\rm Z}}}
\def\IC{\ {\rm I\kern -6.0pt{\rm C}}}

\def\sqr#1#2{{\vcenter{\vbox{\hrule height.#2pt
\hbox{\vrule width.#2pt height #1pt \kern#1pt
\vrule width.#2pt}
\hrule height.#2pt}}}}
\def\square{\mathchoice\sqr56\sqr56\sqr{2.1}3\sqr{1.5}3}

\def\pmb#1{\setbox0=\hbox{$#1$}%
\kern-.025em\copy0\kern-\wd0
\kern.05em\copy0\kern-\wd0
\kern-.025em\raise.0433em\box0 }

\def\pmbb#1{\setbox0=\hbox{$\scriptstyle#1$}%
\kern-.025em\copy0\kern-\wd0
\kern.05em\copy0\kern-\wd0
\kern-.025em\raise.0433em\box0 }

\headline{\hbox to \hsize{\hfil\eightit EACEHL-Jun 9, 1997}}
\vglue1.5truein
\centerline{\bf A MINKOWSKI TYPE TRACE INEQUALITY AND} 
\centerline{\bf STRONG SUBADDITIVITY OF QUANTUM ENTROPY}
\bigskip

{\baselineskip = 12pt
\halign{\qquad#\hfil\qquad\qquad\hfil&#\hfil\cr
Eric  A. Carlen\footnote{$^*$}{Work supported by U.S.
National Science Foundation grant no. DMS 923097} & Elliott H.
Lieb\footnote{$^{**}$}{Work supported by U.S.
National Science Foundation grant no. PHY95--19433--A01.}\cr
School of Mathematics & Departments of Mathematics and Physics\cr
Georgia Institute of Technology& Princeton University\cr
Atlanta, Georgia 30332 & Princeton, New Jersey  08544--0708\cr}
\footnote{$^{}$}{ 
Copyright \copyright 1997 in image and content by the authors. Reproduction of this
article in its entirety by any means is permitted.}}

\vskip 1.2 true in
\centerline{\it Dedicated to M.S. Birman on his 70th birthday.}
\centerline{\it This paper appeared in the  volume in his honor:}
\centerline{\it
Amer. Math. Soc.  Transl. (2), {\rm 189}, 59-69 (1999)}
\bigskip

{\bf Abstract:}  We consider the following trace function on $n$-tuples of
positive operators:
$$\Phi_p(A_1,A_2,\dots,A_n) = Tr\biggl((\sum_{j=1}^n A_j^p)^{1/p}\biggr)$$
and prove that it is jointly concave for $0<p\le 1$ and convex for $p=2$.
We then derive from this a Minkowski type inequality  for operators on a
tensor product of three Hilbert spaces, and show how this implies the
strong subadditivity of quantum mechanical entropy. For $p>2$, $\Phi_p$
is neither convex nor concave. We conjecture that
$\Phi_p$ is convex for $1<p<2$, but our methods do not show this.

\vfill\eject

\def\d{{\rm d}}

\def\h{{\cal H}}
\def\p{{\cal P}}

\noindent{\bf I. INTRODUCTION}
\vskip .3 true in

Let  $\p_\h$ denote the set of all positive semidefinite operators
on a finite dimensional Hilbert space $\h$ with inner product $\langle\cdot,
\cdot\rangle$. 
Then, for any finite natural number $n$, any finite $p>0$, and any finite
$n$-tuple $(A_1,A_2,\dots,A_n)$ of elements of $\p_\h$, define
$$\Phi_p(A_1,A_2,\dots,A_n) = Tr\biggl((\sum_{j=1}^n A_j^p)^{1/p}\biggr)\eqno(1.1)$$

The main result of this paper is the following:
\bigskip
\noindent{\bf Theorem 1}{\it \quad For $0<p\le1$, $\Phi_p$ is a jointly concave function of its 
arguments. For $p=2$, $\Phi_p$ is jointly convex. For $p>2$, $\Phi_p$ is neither
convex nor concave.}
\bigskip

We conjecture that $\Phi_p$ is jointly convex for $1<p<2$. We state all of the theorems in a finite dimensional context, and some of our methods of proof
explicitly involve this finite dimension. Nonetheless, the results themselves
do not depend on the dimension, and therefore easily extend to the appropriate trace classes
on an infinite dimensional Hilbert space.

We note that the trace in Theorem 1 is essential; the asserted trace inequalities {\it do not} hold as operators inequalities. If they did, we would have, for example at $p=2$ that $(A^2+B^2)^{1/2} \le A+B$. This is
of course not true in general for positive operators, as is well known and easily checked.

We shall use Theorem 1 to derive a Minkowski type inequality for traces of
operators on
a product of three Hilbert spaces. To set this in perspective, 
recall that the Minkowski inequality says that for non-negative measurable
functions $f$ on the Cartesian product of two measure spaces $(X,\mu)$ and $(Y,\nu)$,
$$\biggl(\int_X\biggl(\int_Y f(x,y)\d \nu\biggr)^p\d \mu\biggr)^{1/p} \le
\int_Y\biggl(\int_X f^p(x,y)\d \mu\biggr)^{1/p}\d \nu\eqno(1.2)$$
for $p\ge 1$, and that the opposite inequality holds for $0<p\le 1$.

A direct analog of (1.2) holds for positive operators $A$ on the tensor
 product of two 
Hilbert spaces $\h_1\otimes\h_2$. To state it, let
$Tr_1 A$ denote the positive operator on $\h_2$ that is given as a quadratic form by
$$\langle v,Tr_1 A v\rangle = \sum_j\langle u_j\otimes v, A(u_j\otimes v)\rangle$$
where $v\in \h_2$ and the $u_j$ constitute an orthonormal basis of $\h_1$. As is well known, the quadratic form on the left is independent of the choice of the orthonormal basis on the right. The operator
$Tr_1 A$ so defined is called the {\it partial trace of $A$ over $\h_1$}. It will be
convenient, and generally clearer, in what follows to write $Tr_1$ also to denote the usual trace on $\h_1$ for operators $A$ on $\h_1$ alone.

The following is the tracial analog of (1.2):
\bigskip
\noindent{\bf Theorem 2}{\it\quad Let $A$  be a positive operator on the tensor
product of two 
Hilbert spaces $\h_1\otimes\h_2$. Then for all $p\ge 1$,
$$\bigl(Tr_2(Tr_1 A)^p\bigr)^{1/p} \le Tr_1\bigl(\bigl(Tr_ 2 A^p\bigr)^{1/p}\bigr)
\eqno(1.3)$$
and inequality (1.3) reverses for $0<p\le 1$.}

Returning to (1.2), note that it has a trivial extension to functions of three (or more)
variables. Though trivial, it has an interesting consequence.
If one considers a non-negative measurable function $f(x,y,z)$ on the Cartesian
product of three measure spaces $(X,\mu)$, $(Y\nu)$ and $(Z,\rho)$, and simply
holds $z$ fixed as a parameter, one gets
$$\biggl(\int_X\biggl(\int_Y f(x,y,z)\d \nu\biggr)^p\d \mu\biggr)^{1/p} \le
\int_Y\biggl(\int_X f^p(x,y,z)\d \mu\biggr)^{1/p}\d \nu\eqno(1.4)$$
{\it pointwise in $z$} for $p\ge 1$.
Integrating in $z$ then yields
$$\int_Z\biggl(\int_X\biggl(\int_Y f(x,y,z)\d \nu\biggr)^p\d \mu\biggr)^{1/p}\d \rho \le
\int_Z\int_Y\biggl(\int_X f^p(x,y,z)\d \mu\biggr)^{1/p}\d \nu\d \rho\eqno(1.5)$$
for $p\ge 1$, and of course the inequality reverses for $0<p\le1$. 

Now, since (1.5) is an equality at $p=1$, we get another inequality by
differentiating (1.5) with respect to $p$ at $p=1$. This yields an entropy inequality. In fact, using the homogeneity of (1.5), we can normalize $f$ so that it is
a probability density. Recall that for any probability
density $\rho$ on any measure space $(X,\mu)$, 
the entropy $S(\rho)$ is defined as
$$S(\rho) = -\int_X \rho\ln \rho\d \mu\ .\eqno(1.6)$$
We then denote various marginal densities of $f$ as follows:
$$f_{2,3}(y,z) = \int_X f(x,y,z)\d \mu\qquad f_{1,3}(x,z) = \int_Y f(x,y,z)\d \nu \qquad
f_{3}(z) = \int_X\int_Y f(x,y,z)\d \mu\d \nu$$
Then the derivative of (1.5) at $p=1$ is
$$S(f_{1,3}) + S(f_{2,3}) \ge S(f_{1,2,3}) + S(f_{3})\eqno(1.7)$$
which is the strong subadditivity of the classical entropy; see [L75].

Now consider operators on the product of three Hilbert spaces, and a density
matrix $A$; i.e., a positive operator on $\h_1\otimes \h_2\otimes \h_3$
with $Tr A = 1$. The entropy $S(A)$ of a density matrix $A$ is defined by
$$S(A) = -Tr(A\ln A)\ .\eqno(1.8)$$

The operator analog of (1.7) is the Lieb--Ruskai [LR] strong subadditivity 
inequality for the quantum mechanical entropy:

$$S(A_{1,3}) + S(A_{2,3})  \ge S(A_{1,2,3}) + S(A_3)$$
where, in analogy with our notational conventions for marginal densities, we
define
$$A_{1,2,3} = A\,\quad A_{2,3} = Tr_1A\ ,\quad A_3 = Tr_1Tr_2 A$$
and so forth.

Thus, the differential form of Minkowski type inequality (1.7) is known to hold at $p=1$ for operators.
It is therefore natural to enquire whether there exists an operator analog of
the three--variable Minkowski inequality (1.7) for other values of $p$. Unfortunately, the methods at our
disposal suffice to establish this only for $0<p\le 1$ and for $p=2$. 

\bigskip
\noindent{\bf Theorem 3}{\it \quad Let $A$  be a positive operator on the tensor
product of three 
Hilbert spaces $\h_1\otimes\h_2\otimes\h_3$. Then
$$Tr_3\bigl(Tr_2(Tr_1 A)^p\bigr)^{1/p} \le Tr_{1,3}\bigl(\bigl(Tr_ 2 A^p\bigr)^{1/p}\bigr)
\eqno(1.10)$$
for $p=2$ and, trivially, $p=1$, while the reverse inequality holds for $0<p\le 1$.}

\bigskip
This is,
nonetheless, enough to imply the strong subadditivity (1.9): one simply
takes the left derivative at $p=1$. 

It is readily seen by considering block--diagonal matrices that the inequality of Theorem 3 implies the convexity of $\Phi_p$ for $p=2$, and the concavity of $\Phi_p$ for  $0<p\le 1$.
By the same token, (1.10) cannot hold in general for $p>2$ since 
this would imply the convexity of $\Phi_p$ for
such $p$, and Theorem 1 precludes this. This is in contrast to Theorem 2, the Minkowski inequality for
two spaces, which holds for all $p\ge 1$.

The fact that there is such an easy passage from the Minkowski inequality 
in two variables to that in three variables may leave one surprised that there
should be any difficulty in making the same passage with operators. But difficulty there is.
In fact, even the simple version in Theorem 2 seems to require a more 
intricate proof than does 
the corresponding statement for integrals -- which after all
is simply the statement that the unit ball in $L^p$ is convex for $p\ge 1$.
In fact, we know of no previous proof of Theorem 2.

We emphasize that there is no operator analog of the pointwise inequality (1.4).
That is, if we omit $Tr_3$ on both sides of (1.10),  the result will be two operators
on $\h_3$, and these two operators do not satisfy the corresponding operator inequality.

We present a proof of Theorem 2 in Section II. Then in Section III
we prove Theorem 1. In Section IV, we recast Theorem 1 into an equivalent form, from
which Theorem 3 is readily derived in Section IV. Section V contains a brief
comment on a relation between the conjectured convexity for $1<p<2$ and a very interesting
trace inequality of Birman, Koplienko and Solomyak [BKS].

\bigskip
\noindent{\bf II. Proof of Theorem 2}
\bigskip

The following proof of Theorem 2 is given for matrices, but is easily extended to
operators as the statement is dimension independent.

\bigskip

Let $A$ be a positive operator on 
$\p_{\h_1\otimes\h_2}$, the tensor product of two finite dimensional Hilbert spaces.
Suppose first that $p>1$. We proceed by duality.
There is a positive operator $B$ in $\p_{\h_2}$
with $(Tr_2(B^q)^{1/q}) = 1$ with $1/q+1/p = 1$ such that
$$\eqalign{&\bigl(Tr_2(Tr_1 A)^p\bigr)^{1/p} = Tr_2\bigl(B Tr_1A\bigr) = Tr_{1,2}\bigl((I\otimes B)A\bigr) =\cr
&\sum_{i,j}\langle u_i\otimes v_j ,(I\otimes B)A u_i\otimes v_j\rangle =
\sum_{i,j}\langle u_i\otimes B v_j ,A u_i\otimes v_j\rangle\cr}$$
for any pair of orthonormal bases $\{u_i\}$ and $\{v_j\}$. We now choose
the $\{v_j\}$ to be a basis of eigenvectors of $B$, and
let $\{\lambda_j\}$ be the corresponding eigenvalues. 

Then the right hand side above
becomes
$$\eqalign{&\sum_{i,j}\lambda_j\langle u_i\otimes v_j ,A (u_i\otimes v_j)\rangle\le\cr
\biggl(&\sum_j\lambda_j^{q}\biggr)^{1/q}\sum_i\biggl(\sum_j\biggl(\langle
u_i\otimes v_j ,A (u_i\otimes v_j) \rangle\biggr)^p\biggr)^{1/p} =\cr
&\sum_i\biggl(\sum_j\biggl(\langle
u_i\otimes v_j ,A (u_i\otimes v_j) \rangle\biggr)^p\biggr)^{1/p}\cr}$$

Next, by the spectral theorem, for each $i$ and $j$,
$$\langle u_i\otimes v_j ,A (u_i\otimes v_j)\rangle \le
\bigl(\langle u_i\otimes v_j ,A^p (u_i\otimes v_j)\rangle\bigr)^{1/p}$$
Using this, one arrives at
$$\eqalign{\bigl(Tr_2(Tr_1 A)^p\bigr)^{1/p} \le
&\sum_i\biggl(\sum_j\biggl(\langle
u_i\otimes v_j ,A^p(u_i\otimes v_j) \rangle\biggr)\biggr)^{1/p} =\cr
&\sum_i\bigl(\langle
u_i,Tr_2 A^p u_i \rangle\bigr)^{1/p}\cr}$$
Now we choose the $\{u_i\}$ to be a basis of eigenvectors of $Tr_2 A^p$. Then
$$\sum_i\bigl(\langle
u_i,Tr_2 A^p u_i \rangle\bigr)^{1/p} = \sum_i\langle u_i\bigl(Tr_2 A^p\bigr)^{1/p}u_i\rangle
= Tr\bigl(Tr_2 A^p\bigr)^{1/p}$$
and the desired inequality is proved for $p\ge 1$.
Note that this part of the proof works for all $p\ge1$, not only $1\le p \le 2$.

Now suppose $0<p\le 1$, and define $r=1/p$ and $B = A^p$ so that 
$A = B^r$. Since $r>1$, the inequality proved above says
$Tr_1\bigl(\bigl(Tr_2B^r\bigr)^{1/r}\bigr) \ge \bigl(Tr_2\bigl(Tr_1B\bigr)^r\bigr)^{1/r}$.
Rewriting this in terms of $A$ and $p$, and switching the roles
of $\h_1$ and $\h_2$, one obtains the desired result for $0<p\le 1$.\quad $\square$
\bigskip

\noindent{\bf III. Proof of Theorem 1}\qquad 
\bigskip
As before we give the proof for matrices. Consider first the case $0<p<1$. The proof
in this case proceeds by reduction to a theorem of Epstein [E] concerning the function
$$A\mapsto Tr\bigl((BA^pB)^{1/p}\bigr)$$
on $\p_\h$ where $B$ is any given element of $\p_\h$. Epstein's theorem
says that this function is concave for $0<p<1$.

To apply this, consider first the case $n=2$ in (1.1), and define
$${\cal A} =  \left[\matrix{A_1&0\cr 0&A_2\cr}\right]$$
and
$$\sigma =  \left[\matrix{0&I\cr I&0\cr}\right]$$
Then
$${\cal A}^p + \sigma {\cal A}^p\sigma = 
\left[\matrix{A_1^p+A_2^p&0\cr 0&A_1^p+A_2^p\cr}\right]$$
But
$${\cal A}^p + \sigma {\cal A}^p\sigma =
2\biggl({I+\sigma\over 2}\biggr){\cal A}^p\biggl({I+\sigma\over 2}\biggr) +
2\biggl({I-\sigma\over 2}\biggr){\cal A}^p\biggl({I-\sigma\over 2}\biggr)$$
Now define 
$$\Pi_\pm = {I\pm \sigma\over 2}$$
and observe that these are complementary orthogonal projections. Thus,
$$2Tr\bigl((A_1^p + A_2^p)^{1/p}\bigr) = 
2^{1/p}Tr\bigl(\bigl(\Pi_+ {\cal A}^p\Pi_+\bigr)^{1/p}\bigr) +
2^{1/p}Tr\bigl(\bigl(\Pi_- {\cal A}^p\Pi_-\bigr)^{1/p}\bigr)\eqno(3.1)$$
Epstein's theorem, with $A = {\cal A}$ and $B = \Pi_\pm$ now implies that 
each term on the right hand side of (3.1) is a concave function
of ${\cal A}$, which means that the left hand side is a jointly concave function
of $A_1$ and $A_2$. This concludes the proof for $n=2$.

One now easily iterates this procedure to obtain the result for all dyadic
powers $n = 2^k$, and hence for all $n$.

To prove the convexity of $\Phi_2$ there are several way to proceed, but the simplest
was pointed out to us by S. Sahi. Namely, let $n$ be given and consider the block matrix ${\cal A}$ given by
$${\cal A} = 
\left[\matrix{A_1&0&\dots&0\cr A_2&0&\dots &0\cr
\vdots&\vdots&\dots&\vdots\cr
A_n&0&\dots & 0\cr}\right]$$

Then 
$$\Phi_2(A_1,A_2,\dots,A_n) = Tr|{\cal A}|$$
where $|X|$ is the usual operator absolute value; i.e., $\sqrt{X^*X}$.
In other words $\Phi_2(A_1,A_2,\dots,A_n)$ is simply the trace norm of ${\cal A}$,
is therefore clearly jointly convex in $A_1,A_2, \dots,A_n$.

Finally, we show that convexity fails to hold for $p>2$. 
To see this, choose any pair $A_1$,$A_2\in \p_\h$, and any vector $v$ such that
$$\langle v,((A_1^p + A_2^p)/2) v\rangle 
< \langle v,((A_1 + A_2)/2)^p v\rangle\ .\eqno(3.2)$$
Note the strict inequality here. It is always possible to find such $A_1$,$A_2$ and $v$ for $p>2$ since, for such $p$, $X\mapsto X^p$ is not operator convex.

Now let $\Pi_v$ denote the orthogonal projection onto the span of $v$, and let
$\Pi_v^\perp = I - \Pi_v$ denote its orthogonal complement. Then, for a large number $\lambda$ to
be fixed below, put
$$B = \Pi_v + \lambda \Pi_v^\perp\ .$$ Then, if $\Phi_p$ were convex, we would have
$$\limsup_{t\to 0}pt^{-p}\biggl(\Phi_p\biggl(t{A_1+A_2\over 2},B\biggr) -
{1\over 2}\Phi_p(tA_1,B) - {1\over 2}\Phi_p(tA_2,B)\biggr) \le 0\ .\eqno(3.3)$$
However, for small $t>0$,
$$\eqalign{&\Phi_p\biggl(t{A_1+A_2\over 2},B\biggr) = Tr\biggl(t^p\biggl({A_1+A_2\over 2}\biggr)^p + B^p\biggr)^{1/p} =\cr
&TrB + {t^p\over p}Tr\biggl(B^{1-p}\biggl({A_1+A_2\over 2}\biggr)^p\biggr) + O(t^{2p})\cr}$$
and
$$\eqalign{&{1\over 2}\Phi_p(tA_1,B) + {1\over 2}\Phi_p(tA_2,B) =\cr
&TrB + {t^p\over p}\biggl({1\over 2}TrB^{1-p}A_1^p + {1\over 2}TrB^{1-p}A_2^p\biggr) + O(t^{2p})\ .\cr}$$
Thus,
$$\eqalign{&\limsup_{t\to 0}pt^{-p}\biggl(\Phi_p(t{A_1+A_2\over 2},B) -
{1\over 2}\Phi(tA_1,B) - {1\over 2}\Phi(tA_2,B)\biggr) =\cr
&Tr\biggl(B^{1-p}\biggl({A_1+A_2\over 2}\biggr)^p\biggr) -
\biggl({1\over 2}TrB^{1-p}A_1^p + {1\over 2}TrB^{1-p}A_2^p\biggr) =\cr
&\langle v,\biggl({A_1 + A_2\over 2}\biggr)^p v\rangle -
\biggl({\langle v, A_1^pv\rangle\over 2} + {\langle v, A_2^pv\rangle\over 2}\biggr) + O(\lambda^{1-p})\cr}$$
Now taking $\lambda$ sufficiently large, this last term on the right is
stricly positive by (3.2).
This contradicts (3.3), and thus convexity does not hold -- not even separately.
\quad $\square$
\bigskip

\noindent{\bf IV. Corollary of Theorem 1 and Proof of Theorem 3}\qquad 
\bigskip

A corollary of Theorem 1 is obtained by writing the partial trace as an average, and
exploiting the convexity and concavity established above.
Let $A$ be a positive
operator on $\h_1\otimes\h_2$. Next, suppose the dimension of
$\h_2$ is $N$, and fix some orthonormal basis $\{e_1,e_2,\dots,e_N\}$. With respect to this basis,
define the self-adjoint unitary operators $U_{i,j}$ and $V_{i,j}$
on $\h_2$ by
$$\eqalign{&U_{i,j} = I - E_{i,i} - E_{j,j} + E_{i,j} + E_{j,i}\cr
&V_{i} = I - 2E_{i,i} \cr}$$ 
where the $i$ and $j$
are a distinct pair of indices. Let ${\cal G}$ be the subgroup of the group of unitary operators on $\h_2$ that is generated by this family together with the
identity. Each operator $W$ in this group acts by 
$$We_j = (-1)^{s(j)}e_{\pi(j)}$$
for some permutation $\pi(\cdot)$, and some map $s:{1,2,\dots,N}\mapsto {0,1}$.
Thus, the size of the group is $2^NN!$, and the point about it is that any operator on $\h_2$ that commutes with every element of this group is necessarily
a multiple of the identity on $\h_2$.
Then
$${1\over 2^NN!}\sum_{W\in{\cal G}} (I\otimes W^*)A(I\otimes W) = {1\over N}
Tr_2(A)\otimes I_{\h_2}\eqno(4.1)$$
This way of writing partial traces can be traced back to Uhlmann [U]. {}From here one easily arrives at the following result:
\bigskip

\noindent{\bf Theorem 4}\quad{\it For $p>0$,
let the map $\Psi_p(A)$ from positive operators ${\cal A}$ on $\h_1\otimes\h_2$ 
to $\IR_+$ be given by
$$\Psi_p(A) = Tr_1\bigl(\bigl(Tr_2A^p\bigr)^{1/p}\bigr)\ .\eqno(4.2)$$
Then this map is concave for $0<p\le 1$, convex for $p=2$, and neither for
$p>2$.}
\bigskip

\noindent{\bf Proof:}\quad We shall assume that the dimension of $\h_2$ is
$N$ so that we may apply the averaging formula introduced above. We then have
$$\eqalign{&Tr_1\bigl(\bigl(Tr_2A^p\bigr)^{1/p}\bigr) = N^{1/p-1}
Tr_{1,2}\biggl(\biggl({1\over N}Tr_2A^p\otimes I_{\h_2}\biggr)^{1/p}\biggr) =\cr
&N^{1/p-1}Tr_{1,2}\biggl(\biggl(   
{1\over 2^NN!}\sum_{W\in{\cal G}} (I\otimes W^*)A^p(I\otimes W)  \biggr)^{1/p}\biggr)
=\cr
&N^{1/p-1}\biggl({1\over 2^NN!}\biggr)^{1/p}Tr_{1,2}\biggl(\biggl(   
\sum_{W\in{\cal G}}\bigl( (I\otimes W^*)A(I\otimes W) 
\bigr)^p \biggr)^{1/p}\biggr)\cr}
$$
The result now follows directly from Theorem 1.\quad $\square$
\bigskip

Notice that the conclusion of Theorem 4 not only follows from Theorem 1, but also implies it. To see this, suppose that the
$A$ in Theorem 4 is block diagonal with
$$A = 
\left[\matrix{A_1&0&\dots&0\cr 0&A_2&\dots &0\cr
\vdots&\vdots&\dots&\vdots\cr
0&0&\dots & A_n\cr}\right]$$
Then clearly 
$$\Psi_p(A) = \Phi_p(A_1,A_2,\dots,A_n)$$

We remark that if our convexity conjecture turns out to be true for $1<p<2$,
then a proof along the same lines as the proof above of Theorem 4 will prove
the conjectured tracial generalization of Minkowski's inequality.
\bigskip

\noindent{\bf Proof of Theorem 3}\quad Suppose that the dimension of $\h_1$ is $N$. The left hand side of (1.10) can be
written in terms of $\Psi_p$, namely
$$Tr_3\bigl(Tr_2(Tr_1A)^p\bigr)^{1/p} = Tr_{1,3}\biggl(Tr_2\biggl({1\over N}Tr_1A\otimes I_{\h_1}\biggr)^p\biggr)^{1/p} =$$
$$\Psi_p\biggl({1\over N}Tr_1 A\otimes I_{\h_1}\biggr)$$
where the pair of spaces in the definition of $\Psi_p$ is taken to be
$\h_2$ and $\h_1\otimes\h_3$.

Then by (4.1) and the convexity of $\Psi_p$ established in Theorem 4,
$$\eqalign{\Psi_p\biggl({1\over N}Tr_1 A\otimes I_{\h_1}\biggr) = 
&\Psi_p\biggl(   
{1\over 2^NN!}\sum_{W\in{\cal G}} (I\otimes W^*)A(I\otimes W)  \biggr) \le\cr
{1\over 2^NN!}\sum_{W\in{\cal G}}&\Psi_p\biggl((I\otimes W^*)
A(I\otimes W)  \biggr)\cr}$$
The last term above is
$${1\over 2^NN!}\sum_{W\in{\cal G}}Tr_{1,3}\biggl(Tr_2\bigl((I\otimes W^*) A^p (I\otimes W)\bigr)\biggr)^{1/p} = $$
$${1\over 2^NN!}\sum_{W\in{\cal G}}Tr_{1,3}\biggl((I\otimes W^*)\bigl(Tr_2 A^p\bigr) (I\otimes W)\biggr)^{1/p} = $$
$${1\over 2^NN!}\sum_{W\in{\cal G}}Tr_{1,3}\biggl((I\otimes W^*)\bigl(Tr_2 A^p\bigr)^{1/p} (I\otimes W)\biggr) = $$
$$Tr_3 Tr_1\bigl(\bigl(Tr_2 A^p\bigr)^{1/p}\bigr)\ ,$$
which is the desired result.\quad$\square$ 
\bigskip
\noindent{\bf V. The BKS Inequality and the $1<p<2$ Conjecture}
\bigskip

Birman, Koplienko and Solomyak [BKS] proved that for $p>1$, and $A$ and $B$
positive semidefinite operators,
$$Tr\bigl(B^p - A^p\bigr)^{1/p}_+ \ge Tr(B-A)_+\eqno(5.1)$$
where $X_+$ denotes the positive part of a self adjoint operator $X$; i.e., $X_+ = (X+|X|)/2$. In
(5.1), neither $B$ nor $A$  needs to be bounded, but it is assumed that 
$\bigl(B^p - A^p\bigr)^{1/p}_+$ is trace class. Though the inequality in (5.1) is only one of several very interesting inqualities proved in [BKS], 
we refer to it here as the BKS inequality.

The proof is in two parts, the first of which is to reduce consideration to the case
$B^p\ge A^p$ in which case one has $B = (A^p + C^p)^{1/p}$ with $C>0$.
Then (5.1) becomes
$$Tr(C+A) \ge Tr(A^p + C^p)^{1/p}\eqno(5.2)$$
for all $A\ge 0$ and $C\ge 0$. It is (5.2) that interests us here.

Clearly (5.2) can be rewritten as
$$\Phi_p(A,C) \le \Phi_p(A,0) + \Phi_p(0,C)\ ,\eqno(5.3)$$
which is a subadditivity property of $\Phi_p$ for all $p>1$. Since  $\Phi_p$
is homogeneous of degree 1, subadditivity and convexity are the same thing.
Thus for $p=2$, (5.3) is a special case of the convexity of $\Phi_p$ proved in Theorem 1, and for $1<p<2$, it would be a consequence of the conjectured convexity
for these $p$. However, the BKS inequality holds for all $p>1$, not only for
$1<p<2$.

There is a simple proof of (5.2) for matrices. Let
$$
M_\pm = 
\left[\matrix{A^{p/2}&0\cr \pm C^{p/2}&0\cr
}\right]$$
so that
$$Tr\bigl(M_\pm^*M_\pm\bigr)^{1/p} = Tr(A^p+C^p)^{1/p}\ . $$
On the other hand, the spectrum of $M_\pm^* M_\pm$ is the same as
the spectrum of $M_\pm M_\pm^*$, so
$$Tr(A+C) = Tr\bigl(M_\pm M_\pm^*\bigr)^{1/p}\ . $$
One computes
$$M_\pm M_\pm^* =
\left[\matrix{A^p&\pm J\cr \pm J &C^p\cr
}\right]$$
with $J = A^{p/2}C^{p/2}$. Since $X\mapsto Tr(X^{1/p})$ is concave for $p>1$, one has that
$$
Tr(A+C) = Tr\biggl({M_+ M_+^* + M_- M_-^* \over 2}\biggl)^{1/p}
\ge Tr\bigl(A^p + C^p\bigl)^{1/p}\ .$$

A recent application of the BKS inequality, and a different proof of (5.1) that holds in the case of unbounded operators, can be found in [LSS].
\bigskip
\noindent{\bf Acknowledgements}
\bigskip
We thank T. Ando and F. Hiai for a careful reading of this paper, and for pointing out many misprints in an earlier draft.
\bigskip
\noindent{\bf References}
\bigskip

\item{[BKS]} Birman. M.S., Koplienko. L.S. and Solomyak, M.Z.: Estimates for
the spectrum of the difference between fractional powers of two self
adjoint operators, Soviet Mathematics, {\bf 19}(3) 1-6, (1975)

\item{[E]} Epstein, H.: On a concavity theorem of Lieb, Commun. Math. Phys. {\bf 31}
317-327 (1973)

\item{[L73]} Lieb, E.H.: Convex trace functions and the Wigner-Yanase-Dyson Conjecture, Adv. Math. {\bf 11} 267-288 (1973) 

\item{[L75]} Lieb, E.H.: Some Convexity and Subadditivity Properties of Entropy,
Bull. Amer. Math. Soc. {\bf 81} 1-13 (1975) 

\item{[LR]} Lieb, E.H. and Ruskai, M.B.: Proof of the Strong Subadditivity of Quantum-Mechanical Entropy, J. Math. Phys. {\bf 14} 1938-1941 (1973) 

\item{[LSS]} Lieb, E.H., Siedentop, H. and Solovej, J.P.: Stability and
Instability of Relativistic Electrons in Magnetic Fields, J. Stat. Phys.,
(in press).

\item{[U]} Uhlmann, A.: S\"atze \"uber Dichtematrizen, Wiss. Z. Karl-Marx Univ. Leipzig {\bf 20}, 633-53 (1971)

\end